\documentclass[titlepage,a4paper]{article}
\usepackage[T1]{fontenc}
\usepackage[latin1]{inputenc}
\usepackage{amsmath}
\usepackage{amssymb}
\usepackage{amsthm}

\newtheorem{Def}{Definition}[section]

\newtheorem{Th}{Theorem}[section]

\newtheorem{lem}{Lemma}[section]

\newcommand{\RR}{\mathbb{R}}
\newcommand{\EE}{\mathbb{E}}
\newcommand{\PP}{\mathbb{P}}

\newcommand{\QQ}{\mathbb{Q}}

\newcommand{\cF}{{\cal F}}

\newcommand{\cG}{{\cal G}}

\newcommand{{{\cadlag}}}{c\`adl\`ag}

\renewcommand{\baselinestretch}{2}

\begin{document}

\sffamily
\pagestyle{plain}

\title{On measures of unfairness and an optimal currency transaction tax }
\date{}
\author{Frederik S Herzberg\thanks{Mathematical Institute, University of Oxford, Oxford OX1 3LB, United Kingdom, Telephone (reception): +44 (0)1865 273525, Facsimile: +44 (0)1865 273583, e-mail: {\tt herzberg@maths.ox.ac.uk}. The author gratefully acknowledges funding from the German National Academic Foundation ({\em Studienstiftung des deutschen Volkes}) and a PhD grant of the German Academic Exchange Service ({\em Doktorandenstipendium des Deutschen Akademischen Austauschdienstes}).}}

\renewcommand{\baselinestretch}{1} \maketitle

\begin{abstract} The purpose of the present paper is the analysis of a model describing how herd behaviour and self-fulfilling prophecies can influence currency exchange rates, and what the impact of a currency transaction tax would be. These considerations yield a stochastic differential equation that can be studied using the methods of infinitesimal analysis. We will show, using a suitable notion of unfairness for discounted price processes which measures the distance from being a martingale, that the fairest tax rate is the maximal one subject to the condition that it does not affect real-economic speculation.
\end{abstract}

\renewcommand{\baselinestretch}{2}

\section{Introduction}

\renewcommand{\baselinestretch}{2}

The speculative bubbles of the 1990's have led to the question whether there is an intrinsic tendency to herd behaviour on the financial markets. Corcos et al \cite{CEMMS} have analysed these phenomena from a microeconomic point of view, modelling the behaviour of ``bullish'' and ``bearish'' traders; their model is time-discrete, but admits a limit in which the proportion of bullish traders is deterministic. As a very recent study in the general theory of herding we refer to Drehmann, Oechssler and Roider \cite{DOR}.

Our model for the behaviour of the pre-crash stock price will be time-continuous. However, we do not lose the explanatory advantages of a discrete model, since there is a straightforward hyperfinite analogue to the stochastic differential equation that we assume as a model, and the standard part of this hyperfinite object solves the SDE. It is a fairly basic model which assumes exponential growth of the bubbles. Whilst the model devised by Corcos et al \cite{CEMMS} entails exponential growth of the bubbles at the beginning, their model predicts that bullish behaviour prompts singularities which will be muted to a chaotic behaviour due to bearish agents. In our model, we assume that once the bubble has been recognised as not being justifiable in terms of real-economic data, an exponential decrease in the stock price is going to occur.

To the knowledge of the author, no attempt has been made so far to develop quantitative notions of unfairness for price processes. In a forthcoming paper, suitable measures of unfairness will be investigated in more detail.

Our proofs will employ techniques from nonstandard analysis; these have already been successfully applied to mathematical finance in the work of Cutland, Kopp and Willinger \cite{CKW93}.

\section{The basic model} 

As was pointed out in the introduction, the pre-crash stock-price process shall be considered. For this purpose, $\alpha:\RR\rightarrow \RR$ will be a piecewise constant function such that \begin{eqnarray*} \alpha=\alpha(1)>0 \text{ on } (0,+\infty), \quad &&\alpha=\alpha(-1)<0 \text{ on } (-\infty,0).\end{eqnarray*} We assume that the logarithmic discounted price process $x^{(\upsilon)}$ is governed, given some initial condition, by the stochastic differential equation \begin{eqnarray}\label{stock}dx^{(\upsilon)}_t &=& \alpha\left(x^{(\upsilon)}_t- \sum_{u\in I} p_u \cdot x^{(\upsilon)}_{(t-u)\vee 0}\right)\chi_{\left\{\left| x^{(\upsilon)}_t-\sum_{u\in I} p_u \cdot x^{(\upsilon)}_{(t-u)\vee 0}\right|\geq\upsilon \right\}}dt \\ && \nonumber +\sigma\cdot db_t - \frac{{\sigma}^2}{2} dt,\end{eqnarray} 
where $r>0$ is the logarithmic discount rate, $b$ is the one-dimensional Wiener process, $(p_u)_{u\in I}$ is a convex combination -- i.e. $I\subset(0,+\infty)$ is finite, $\forall u\in I\quad p_u>0$ and $\sum_{u\in I}p_u=1$. The parameter $\upsilon$ depends very much on the tax rate we assume. If $\rho $ is the logaritmic tax rate and $T$ the expected time during which one will hold the asset (that is the expected duration of the upward or downward tendency of the stock price), one can compute $\upsilon$ as follows: $$\upsilon=T\cdot \rho.$$ This, of course is first of all only a formal equation that -- thanks to the boundedness of $\alpha$ -- can be made rigorous using the theory of stochastic differential equations developed by Hoover and Perkins \cite{HP83} or Albeverio et al \cite{AFHL86} for instance.
We will introduce the following notation: $$\psi^{(\upsilon)}:=\chi_{\left\{|\cdot|\geq\upsilon\right\}}\cdot \alpha.$$

\section{How to minimise ``unfairness''}

\begin{Def} \label{fair} For any stochastic process $y$ on an adapted probability space $\left(\Gamma,\cG,\QQ\right)$, we define the integral $$m\left(y,\Gamma\right):=\int_0^1\EE_\QQ\left|\left.\frac{1}{y_t}\frac{d}{du}\right|_{u=0}\EE\left[y_{t+u}\left.\right|\cG_t\right]\right|dt$$ (if existent) to be the {\em complete-market unfairness} of $y$. Whenever the derivative does not exist on a set of positive $\lambda^1\otimes\QQ$-measure, we set $m\left(y,\Gamma\right)=+\infty$.
\end{Def}

Intuitively, this function measures how often (in terms of time and probability) and how much it will be the case that one may expect to obtain a multiple (or a fraction) of one's portfolio simply by selling or buying the stock under consideration.

More general notions of unfairness will be studied as part of a forthcoming paper. We will confine our attention to the properties of the notion of unfairness just defined and state the following immediate result which justifies the name if we accept the common view that martingales correspond to fair games (alluding to the well-known as well as intuitively plausible equivalence between non-existence of arbitrage and the martingale property of the discounted price process).

\begin{lem} For all semimartingales $z$ on an adapted probability space $\left(\Gamma,\cG,\QQ\right)$, $z$ is a martingale if and only if $m(z,\Gamma)=0$. Moreover, $m(\cdot,\Gamma)$ remains unchanged under multiplication by constants.
\end{lem}

In order to study $m$, it is advisable to resort to nonstandard methods, for the relatively large number of limit processes involved in the defnition of $m$ cannot be expected to ease a standard analysis of $m$'s properties.

First of all we derive a formula that makes it easier to actually compute $m$ in our concrete setting:

\begin{lem} \label{m() formula}If $x^{(\upsilon)}$ satisfies (\ref{stock}) for some $\upsilon>0$ on an adapted probability space $\left(\Gamma,\cG,\QQ\right)$, then the discounted price process $\left(\exp\left(x^{(\upsilon)}_t\right) \ : t\geq 0\right)$ is of finite unfairness. More specifically, \begin{eqnarray*} m\left(\exp\left(x^{(\upsilon)}\right),\Gamma\right) &=& \int_0^1\EE\left[\left|\psi^{(\upsilon)}\left(x_t^{(\upsilon)}-\sum_{i\in I}p_ix_{(t-i)\vee 0}^{(\upsilon)}\right)\right|\right] dt\\ &=& \int_0^1\EE\left[\left|  \left.\frac{d}{du}\right|_{u=0}\EE\left[\left. x_{t+u}^{(\upsilon)}\right|\cG_t\right]+\frac{\sigma^2}{2} \right|\right]  dt\end{eqnarray*}
\end{lem}
\begin{proof} The proof is more or less a formal calculation, provided one is aware of the path-continuity of our process and the fact that the filtrations generated by $b$ and $x^{(\upsilon)}$ are identical. For this implies that, given $t>0$, the value $\psi^{(\upsilon)}\left(x_{t+u}^{(\upsilon)}-\sum_{i\in I}p_ix_{(t+u-i)\vee 0}^{(\upsilon)}\right)(\omega) $ does not change within sufficiently small times $u$ -- almost surely for all those paths $\omega$ where $x_{t}^{(\upsilon)}(\omega)-\sum_{i\in I}p_ix_{(t-i)\vee 0}^{(\upsilon)}(\omega)\notin \{\pm\upsilon\}$, this condition itself being satisfied with probability $1$. Now, using this result and the martingale property of the quotient of the exponential Brownian motion and its exponential bracket, we can deduce that for all $t>0$ almost surely:
\begin{eqnarray*}&&\frac{1}{\exp\left(x_{t}^{(\upsilon)}\right)}\frac{d}{du}\left|_{u=0}
\EE\left[\left.\exp\left(x_{t+u}^{(\upsilon)}\right)\right|\cG_t\right]\right.\\ &=& \frac{d}{du}\left|_{u=0}
\EE\left[\left.\exp\left(x_{t+u}^{(\upsilon)}-x_{t}^{(\upsilon)}\right)\right|\cG_t\right]\right.\\ &=&  \frac{d}{du}\left|_{u=0} \exp\left(\psi^{(\upsilon)}\left(x_t^{(\upsilon)}-\sum_{i\in I}p_ix_{(t-i)\vee 0}^{(\upsilon)}\right) u\right) \right.\\ && \cdot \EE\left[\left.\exp\left(\sigma b_{t+u}-\frac{\sigma^2}{2}(t+u)-\left(\sigma b_{t}-\frac{\sigma^2}{2}t\right)\right)\right|\cG_t\right] \\ &=&\frac{d}{du}\left|_{u=0}  \exp\left(\psi^{(\upsilon)}\left(x_t^{(\upsilon)}-\sum_{i\in I}p_ix_{(t-i)\vee 0}^{(\upsilon)}\right) u\right) \right. \\ && \cdot \EE\left[\left.\exp\left(\sigma b_{t+u}-\frac{\sigma^2}{2}(t+u)\right)\right|\cG_t\right]\exp\left(-\sigma b_{t}+\frac{\sigma^2}{2}t\right) \\ &=&\frac{d}{du}\left|_{u=0} \exp\left(\psi^{(\upsilon)}\left(x_t^{(\upsilon)}-\sum_{i\in I}p_ix_{(t-i)\vee 0}^{(\upsilon)}\right) u\right)\cdot 1\right. \\ &=& \psi^{(\upsilon)}\left(x_t^{(\upsilon)}-\sum_{i\in I}p_ix_{(t-i)\vee 0}^{(\upsilon)}\right) .\end{eqnarray*}

Analogously, one may prove the second equation in the Lemma: For, one readily has almost surely
\begin{eqnarray*}&& \frac{d}{du}\left|_{u=0}\EE\left[\left.x_{t+u}^{(\upsilon)}\right|\cG_t\right] \right.\\ &=&  \frac{d}{du}\left|_{u=0}\left( \EE\left[\left.\psi^{(\upsilon)}\left(x_t^{(\upsilon)}-\sum_{i\in I}p_ix_{(t-i)\vee 0}^{(\upsilon)}\right) u\right|\cG_t\right] -\frac{\sigma^2}{2}u + x_t^{(\upsilon)} \right)\right.\\ &=& \psi^{(\upsilon)}\left(x_t^{(\upsilon)}-\sum_{i\in I}p_ix_{(t-i)\vee 0}^{(\upsilon)}\right) -\frac{\sigma^2}{2}\end{eqnarray*} 
In order to proceed from these pointwise almost sure equations to the assertion of the Theorem, one will apply Lebesgue's Dominated Convergence Theorem, yielding $$ m\left(\exp\left(x^{(\upsilon)}\cdot\right),\Gamma\right)= \int_0^1\EE\left[\left|  \left.\frac{d}{du}\right|_{u=0}\EE\left[\left. x_{t+u}^{(\upsilon)}\right|\cG_t\right] \right|\right]  dt.$$

\end{proof}

For finite hyperfinite adapted probability spaces, an elementary proof for the main Theorem can be contrived:

\begin{lem} \label{finite} For any hyperfinite number $H$ we will let $X^{(\upsilon)}$ for all $\upsilon>0$ denote the solution to the hyperfinite initial value problem \begin{eqnarray} && X^{(\upsilon)}_0 = 0,\nonumber \\
&&\label{stock-internal}\forall {t\in\left\{0,\dots,1-\frac{1}{H!}\right\}}\nonumber \\&& X^{(\upsilon)}_{t+\frac{1}{H!}}-X^{(\upsilon)}_{t}\nonumber \\ &=& \alpha\left(X^{(\upsilon)}_t- \sum_{u\in I} p_u \cdot X^{(\upsilon)}_{(t-u)\vee 0}\right)\chi_{\left\{\left| X^{(\upsilon)}_t-\sum_{u\in I} p_u \cdot X^{(\upsilon)}_{(t-u)\vee 0}\right|\geq\upsilon \right\}}\cdot \frac{1}{H!} \\ && \nonumber +\sigma\cdot \pi_{t+\frac{1}{H!}}\cdot\frac{1}{(2H!)^{1/2}} - \frac{{\sigma}^2}{2} \frac{1}{H!}\end{eqnarray} (where $\pi_{\ell/H!}:\Omega=\{\pm 1\}^{H!}\rightarrow \{\pm 1\}$ is for all hyperfinite $\ell\leq H!$ the projection to the $\ell$-th coordinate) which is just the hyperfinite analogue to (\ref{stock}). Using this notation, and considering a finite (rather than merely hyperfinite) adapted probability space of mesh size $H!$, one has for all $k<H!$, \begin{eqnarray*}&&\sum_{k<{H!}} \EE \left|X^{(\tau)}_{k/H!}-\EE\left[\left. X^{(\tau)}_{\frac{k+1}{H!}}\right|\cF_{k/H!}\right]-\frac{\sigma^2}{2}\frac{1}{H!}\right| \\ &\leq& \sum_{k<{H!}} \EE \left|X^{(\upsilon )}_{k/H!}-\EE\left[\left. X^{(\upsilon )}_{\frac{k+1}{H!}}\right|\cF_{k/H!}\right]-\frac{\sigma^2}{2}\frac{1}{H!}\right| \end{eqnarray*} for all $\tau\geq \upsilon $. As a consequence, $m\left(X^{(\cdot)},\Omega\right)$ is monotonely decreasing for all finite adapted probability spaces $\Omega=\{\pm 1\}^{H!}$.
\end{lem}

By transfer to the nonstandard universe, we will obtain the same result for infinite hyperfinite $H$ as well.

\begin{proof}[Proof of Lemma \ref{finite}] The proof of this Lemma relies on exploiting the assumption that $\alpha$ is piecewise constant, since $\alpha=\alpha(1)\chi_{(0,+\infty)}+\alpha(-1)\chi_{(-\infty,0)}+\alpha(0)\chi_{\{0\}}$ yields
\begin{eqnarray*}&& \forall k< H! \quad \forall \upsilon>0 \\ &&\EE\left[\left.X^{(\upsilon)}_{\frac{k+1}{H!}}\right|\cF_{k/H!}\right]-X^{(\upsilon)}_{\frac{k}{H!}} +\frac{\sigma^2}{2}\frac{1}{H!}\\ &=& \frac{1}{H!}\alpha\left(X^{(\upsilon)}_{k/H!}-\sum_{u\in I} p_u \cdot X^{(\upsilon)}_{\left(\frac{k}{H!}-u\right)\vee 0}\right)\chi_{\left\{\left|X^{(\upsilon)}_{k/H!}-\sum_{u\in I} p_u \cdot X^{(\upsilon)}_{\left(\frac{k}{H!}-u\right)\vee 0} \right|\geq\upsilon \right\}} \\ &=& \frac{1}{H!}\left(\begin{array}{c}\alpha(1) \chi_{\left\{ X^{(\upsilon)}_{k/H!}-\sum_{u\in I} p_u \cdot X^{(\upsilon)}_{\left(\frac{k}{H!}-u\right)\vee 0} \geq\upsilon \right\}} \\+ \alpha(-1) \chi_{\left\{ X^{(\upsilon)}_{k/H!}-\sum_{u\in I} p_u \cdot X^{(\upsilon)}_{\left(\frac{k}{H!}-u\right)\vee 0} \leq -\upsilon \right\}} \end{array}\right),\end{eqnarray*} which immediately follows from the construction of Anderson's random walk \cite{A76} $B_\cdot=\frac{1}{\sqrt{2H!}}\pi_\cdot$ and the recursive difference equation defining the process $X^{(\upsilon)}$. However, this last equation implies
\begin{eqnarray*}&& \forall k< H! \quad \forall \upsilon>0 \\ &&\EE\left|\EE\left[\left.X^{(\upsilon)}_{\frac{k+1}{H!}}\right|\cF_{k/H!}\right]-X^{(\upsilon)}_{\frac{k}{H!}}+ \frac{\sigma^2}{2}\frac{1}{H!}\right| \\ &=& \frac{1}{H!}\EE\left|\begin{array}{c}\alpha(1) \chi_{\left\{ X^{(\upsilon)}_{k/H!}-\sum_{u\in I} p_u \cdot X^{(\upsilon)}_{\left(\frac{k}{H!}-u\right)\vee 0} \geq\upsilon \right\}} \\+ \alpha(-1) \chi_{\left\{ X^{(\upsilon)}_{k/H!}-\sum_{u\in I} p_u \cdot X^{(\upsilon)}_{\left(\frac{k}{H!}-u\right)\vee 0} \leq -\upsilon \right\}} \end{array}\right| \\ &=&\frac{1}{H!} \left(\begin{array}{c}\left|\alpha(1)\right|\PP\left\{ X^{(\upsilon)}_{k/H!}-\sum_{u\in I} p_u \cdot X^{(\upsilon)}_{\left(\frac{k}{H!}-u\right)\vee 0} \geq\upsilon \right\}\\+ \left|\alpha(-1)\right| \PP\left\{ X^{(\upsilon)}_{k/H!}-\sum_{u\in I} p_u \cdot X^{(\upsilon)}_{\left(\frac{k}{H!}-u\right)\vee 0} \leq -\upsilon \right\}\end{array}\right).\end{eqnarray*} Now all that remains to be shown is that $\sum_{k< H!}\PP\left\{ X^{(\upsilon)}_{k/H!}-\sum_{u\in I} p_u \cdot X^{(\upsilon)}_{\left(\frac{k}{H!}-u\right)\vee 0} \geq\upsilon \right\}$ and $\sum_{k< H!}\PP\left\{ X^{(\upsilon)}_{k/H!}-\sum_{u\in I} p_u \cdot X^{(\upsilon)}_{\left(\frac{k}{H!}-u\right)\vee 0} \leq-\upsilon \right\} $
are monotonely decreasing in $\upsilon$. The former statement is a consequence of the following assertion: \begin{eqnarray}&&\forall \upsilon'\leq \upsilon\forall \ell< H!\nonumber\\ &&\sum_{k\leq \ell}\PP\left\{ X^{(\upsilon)}_{k/H!}-\sum_{u\in I} p_u \cdot X^{(\upsilon)}_{\left(\frac{k}{H!}-u\right)\vee 0} \geq\upsilon \right\} \nonumber \\ &\leq& \sum_{k\leq \ell}\PP\left\{ X^{(\upsilon')}_{k/H!}-\sum_{u\in I} p_u \cdot X^{(\upsilon')}_{\left(\frac{k}{H!}-u\right)\vee 0} \geq\upsilon' \right\} \label{counting in probability}\end{eqnarray} One can prove this estimate by considering the minimal $\ell$ such that the pointwise inequality \begin{eqnarray}&&\sum_{k\leq \ell}\chi_{\left\{ X^{(\upsilon)}_{k/H!}-\sum_{u\in I} p_u \cdot X^{(\upsilon)}_{\left(\frac{k}{H!}-u\right)\vee 0} \geq\upsilon \right\}} \nonumber \\ &\leq& \sum_{k\leq \ell}\chi_{\left\{ X^{(\upsilon')}_{k/H!}-\sum_{u\in I} p_u \cdot X^{(\upsilon')}_{\left(\frac{k}{H!}-u\right)\vee 0} \geq\upsilon' \right\}}\label{counting pathwise}\end{eqnarray} fails to hold. Then one has an $\omega\in\Omega$ such that 
\begin{eqnarray} \forall{k< \ell}&&\chi_{\left\{ X^{(\upsilon)}_{k/H!}-\sum_{u\in I} p_u \cdot X^{(\upsilon)}_{\left(\frac{k}{H!}-u\right)\vee 0} \geq\upsilon \right\}}(\omega)\nonumber \\ &=&\chi_{\left\{ X^{(\upsilon')}_{k/H!}-\sum_{u\in I} p_u \cdot X^{(\upsilon')}_{\left(\frac{k}{H!}-u\right)\vee 0} \geq\upsilon' \right\}}(\omega), \label{before ell}\\ 1&=&\chi_{\left\{ X^{(\upsilon)}_{\ell/H!}-\sum_{u\in I} p_u \cdot X^{(\upsilon)}_{\left(\frac{\ell}{H!}-u\right)\vee 0} \geq\upsilon \right\}}(\omega), \label{in ell}\\ 0 &=& \chi_{\left\{ X^{(\upsilon')}_{\ell/H!}-\sum_{u\in I} p_u \cdot X^{(\upsilon')}_{\left(\frac{\ell}{H!}-u\right)\vee 0} \geq\upsilon' \right\}}(\omega)\label{in ell '}. \end{eqnarray} But equation (\ref{before ell}) implies, via the difference equation for $X^{(\cdot)}$ (\ref{stock-internal}), inductively in $k$ the relation $$\forall{k< \ell}\quad X^{\upsilon}_k(\omega)=X^{\upsilon'}_k(\omega).$$ If one combines this with $$1=\chi_{\left\{ X^{(\upsilon)}_{\ell/H!}-\sum_{u\in I} p_u \cdot X^{(\upsilon)}_{\left(\frac{\ell}{H!}-u\right)\vee 0} \geq\upsilon \right\}}(\omega)$$ (which is equation (\ref{in ell})) and $\upsilon\leq \upsilon'$, one can derive -- again via the recursive difference equation (\ref{stock-internal}) -- that $X^{(\upsilon')}_{\ell/H!}\geq X^{(\upsilon)}_{\ell/H!}$ applied in $t=\frac{\ell}{H!}$) as well as the estimates \begin{eqnarray*}\chi_{\left\{ X^{(\upsilon')}_{\ell/H!}-\sum_{u\in I} p_u \cdot X^{(\upsilon')}_{\left(\frac{k}{H!}-u\right)\vee 0} \geq\upsilon' \right\}}(\omega)&=& \chi_{\left\{ X^{(\upsilon')}_{\ell/H!}-\sum_{u\in I} p_u \cdot X^{(\upsilon)}_{\left(\frac{\ell}{H!}-u\right)\vee 0} \geq\upsilon' \right\}}(\omega)\\ &\geq& \chi_{\left\{ X^{(\upsilon)}_{\ell/H!}-\sum_{u\in I} p_u \cdot X^{(\upsilon)}_{\left(\frac{\ell}{H!}-u\right)\vee 0} \geq\upsilon' \right\}}(\omega)\\ &\geq &\chi_{\left\{ X^{(\upsilon)}_{\ell/H!}-\sum_{u\in I} p_u \cdot X^{(\upsilon)}_{\left(\frac{\ell}{H!}-u\right)\vee 0} \geq\upsilon \right\}}(\omega)\\ &=&1.\end{eqnarray*} This contradicts equation (\ref{in ell '}). Hence, the estimate (\ref{counting pathwise}) has been established for all $k<H!$, leading to (\ref{counting in probability}).

Similarly, one can prove \begin{eqnarray*}&&\forall \upsilon'\leq \upsilon\forall \ell< H!\\ &&\sum_{k\leq \ell}\PP\left\{ X^{(\upsilon)}_{k/H!}-\sum_{u\in I} p_u \cdot X^{(\upsilon)}_{\left(\frac{k}{H!}-u\right)\vee 0} \leq-\upsilon \right\} \\ &\leq& \sum_{k\leq \ell}\PP\left\{ X^{(\upsilon')}_{k/H!}-\sum_{u\in I} p_u \cdot X^{(\upsilon')}_{\left(\frac{k}{H!}-u\right)\vee 0} \leq-\upsilon' \right\}\end{eqnarray*} which entails that $\sum_{k< H!}\PP\left\{ X^{(\upsilon)}_{k/H!}-\sum_{u\in I} p_u \cdot X^{(\upsilon)}_{\left(\frac{k}{H!}-u\right)\vee 0} \leq-\upsilon \right\} $ must be monotonely decreasing in $\upsilon$.
\end{proof}

Using nonstandard analysis and the model theory of stochastic processes as developed by Keisler and others \cite{HK84,K88,FK02}, we can prove the following result:

\begin{Th} Suppose $\left(y^{(\upsilon)} \ : \ \upsilon>0\right)$ is a family of stochastic processes on an adapted probability space $\Gamma$ such that $y^{(\upsilon)}$ solves the stochastic differential equation (\ref{stock}) formulated above for all $\upsilon>0$. Then the function $\sigma\mapsto m\left(y^{(\sigma)},\Gamma\right)$ attains its minimum on $[0,S]$ in $S$.
\end{Th}
\begin{proof} By the previous Lemmas \ref{finite} and the formula for $m$ from Lemma \ref{m() formula}, the assertion of the Theorem holds true for finite adapted probability spaces. By transfer to the nonstandard universe, we obtain the same result internally for hyperfinite adapted spaces. Now, according to results by Hoover and Perkins \cite{HP83} as well as Albeverio et al \cite{AFHL86}, the solution $X^{(\upsilon)}$ of the hyperfinite initial value problem (\ref{stock-internal}) is a lifting for the solution $x^{(\upsilon)}$ of (\ref{stock}) on a hyperfinite adapted space for any $\upsilon \geq 0$. Now, $y\mapsto m\left(y,\Omega\right)$ is the expectation of a conditional process in the sense of Fajardo and Keisler \cite{FK02}. Therefore, due to the Adapted Lifting Theorem \cite{FK02}, we must have $$\forall \upsilon\geq 0 \quad °m\left(X^{(\upsilon)},\Omega\right) =m\left(x^{(\upsilon)},\Omega\right)$$ (where we identify $m$ with its internal analogue when applied to internal processes). Since the internal equivalent of the Theorem's assertion holds for internal hyperfinite adapted space, the previous equation implies that it is also true for Loeb hyperfinite adapted spaces. \\
Now let $\left(y^{(\upsilon)} \ : \ \upsilon>0 \right)$ be a family of processes on some (not necessarily hyperfinite) adapted probability space $\Gamma$ with the properties as in the Theorem. Because of the universality of hyperfinite adapted spaces \cite{FK02,K88}, we will find a process $x^{(\upsilon)}$ on any hyperfinite adapted space $\Omega$ such that $x$ and $y$ are automorphic to each other. This implies \cite{AFHL86} that $x^{(\upsilon)}$ satisfies (\ref{stock}) as well. Furthermore, as one can easily see using Lebesgue's Dominated Convergence Theorem, $$\forall \upsilon>0 \quad m\left(x^{(\upsilon)},\Omega\right)=m\left(y^{(\upsilon)},\Gamma\right).$$ Due to our previous remarks on the solutions of (\ref{stock}) on hyperfinite adapted spaces, this suffices to prove the Theorem. 
\end{proof}

\section{Conclusions}

We have introduced a suitable notion of the unfairness which a discounted asset price process might have (whether this process deserves to be called ``price'' even though it is not a martingale is a question of terminology, we are concerned with the actual observed amount of money to be paid at a stock exchange). We also have assumed a model of currency prices where extrapolation from the observed behaviour of other agents is, up to white noise and constant inflation, completely accounts for the development of the currency price. In analogy to a Nash equilibrium, we assume that every agent is acting in such a manner that he gains most if all other agents follow his pattern. In our model this pattern will consist in using some average value of past currency prices as a ``sunspots'', that is a proxy for a general perception that depreciation or appreciation of the particular currency in question is due.

In addition, we have assumed that a transaction tax of constant rate has to be paid for some currency transactions. Needless to say, it is in many cases difficult to determine and establish beyond doubt whether a particular currency transaction was based on rigorous analysis of fundamental economic data or merely on expectations on the market. It is therefore reasonable to demand that there be a moderate upper bound on any currency transaction tax rate such that the transaction tax does not prevent any speculation based on fundamental economic data (as opposed to speculation on the behaviour of fellow traders).

In such a setting we have proven that as long as one does not exceed the said upper bound, any rise in the currency transaction tax enhances fairness of the currency market in the technical sense of Definition \ref{fair}.


\begin{thebibliography}{20}

\bibitem{A86} S Albeverio, {\em Some personal remarks on nonstandard analysis in probability theory and mathematical physics}, (Proceedings of the) VIIIth International Congress on Mathematical Physics -- Marseille 1986 (eds. M. Mebkhout, R. S\'en\'eor), World Scientific, Singapore 1987, 409 -- 420.
\bibitem{AFHL86} S Albeverio, J Fenstad, R H\o egh-Krohn, T Lindstr\o m, {\em Nonstandard methods in stochastic analysis and mathematical physics}, Academic Press, Orlando 1986.
\bibitem{A76} R M Anderson, {\em A nonstandard representation for Brownian motion and It\^o integration}, Israel Journal of Mathematics {\bf 25} (1976), 15 -- 46.
\bibitem{CEMMS} A Corcos et al, {\em Imitation and contrarian behaviour: hyperbolic bubbles, crashes and chaos}, Quantitative Finance {\bf 2} (2002), 264 -- 281.
\bibitem{CKW93} N Cutland, E Kopp, W Willinger, {\em A nonstandard treatment of options driven by Poisson processes}, Stochastics and Stochastics Reports {\bf 42} (1993), 115 -- 133.
\bibitem{DOR} M Drehmann, J Oechssler, A Roider, {\em Herding with and without payoff externalities -- an internet experiment}, mimoe 2004.
\bibitem{FK02} S Fajardo, H J Keisler, {\em Model theory of stochastic processes}, Lecture Notes in Logic 14, A. K. Peters, Natick (Mass.) 2002.
\bibitem{HK84} D N Hoover, H J Keisler, {\em Adapted probability distributions}, Transactions of the American Mathematical Society {\bf 286} (1984), 159 -- 201.
\bibitem{HP83} D N Hoover, E Perkins, {\em Nonstandard construction of the stochastic integral and applications to stochastic differential equations I, II}, Transactions of the American Mathematical Society {\bf 275} (1983), 1 -- 58. 
\bibitem{K88} H J Keisler, {\em Infinitesimals in probability theory}, Nonstandard analysis and its applications (ed N Cutland), London Mathematical Society Student Texts 10, Cambridge University Press, Cambridge 1988, 106 -- 139.
\bibitem{L75} P A Loeb, {\em Conversion from nonstandard to standard measure and applications in probablility theory}, Transactions of the American Mathematical Society {\bf 211} (1974), 113 -- 122.
\bibitem{N87} E Nelson, {\em Radically elementary probability theory}, Annals of Mathematics Studies 117, Princeton University Press, Princeton (NJ) 1987.
\bibitem{P83} E Perkins, {\em Stochastic processes and nonstandard analysis}, Nonstandard analysis -- recent developments (ed A E Hurd), Lecture Notes in Mathematics 983, Springer, Berlin 1983, 162 -- 185.
\bibitem{SB86} K D Stroyan, J M Bayod, {\em Foundations of infinitesimal stochastic analysis}, Studies in Logic and the Foundations of Mathematics 119, North-Holland, Amsterdam 1986.

\end{thebibliography}
\end{document}